\documentstyle[12pt]{amsart}
\newcommand{\g}{\goth}
\newcommand{\gtg}{\mbox{\g g}}

\newcommand{\hgtg}{\mbox{$\hat{\gtg}$}}

\newcommand{\vep}{\varepsilon}

\theoremstyle{plain}
 \newtheorem{thm}{Theorem}[section]

\theoremstyle{definition}
 \newtheorem{defn}{Definition}[section]
 
 \newtheorem{ex}{Example}[section]
\theoremstyle{remark}

 \newtheorem{ack}{Acknowledgment}
      
\begin{document}
\title{Hopf algebra extension of a Zamolochikov algebra 
and its double}
\author{Jintai Ding}
 
\address{Jintai Ding, RIMS, Kyoto University}

\maketitle
\begin{abstract}
The  particles with a scattering matrix $R(x)$ are defined as operators
$\Phi_i(z)$ satisfying the relation 
$ R_{i,j}^{j',i'}(x_1/x_2) \Phi_{i'}(x_1)\Phi_{j'}(x_2)=
\Phi_i(x_2)\Phi_j(x_1)$. The algebra generated by those operators 
is  called a Zamolochikov algebra. We 
construct a new Hopf algebra by adding half of the 
FRTS construction of a quantum affine algebra with 
this $R(x)$. Then we double it to obtain a new Hopf algebra such that 
the full FRTS construction of a quantum affine algebra 
is a Hopf subalgebra inside. 
Drinfeld realization of quantum affine algebras is included as 
an  example. This is a further generalization of the constructions
in \cite{DI}. 
\end{abstract}

\pagestyle{plain}
\section{Introduction.} 

In physics, the particles with a scattering matrix $R(x)$ in 
End($V)\otimes$ End($V)$  are defined
with the operators $\Phi_i(x)$ index by a linear independent basis
of $V$ such that 
$$ R_{i,j}^{j',i'}(x_1/x_2) \Phi_{i'}(x_1)\Phi_{j'}(x_2)=
\Phi_i(x_2)\Phi_j(x_1), $$
where $V$ is a vector space, $x$ is a parameter in $\Bbb C$. 
This naturally gives a algebra with these  current generators
$\Phi_i(x)$, which we will call a  Zamolochikov algebra.  However 
this algebra  is not given 
a Hopf algebra structure. We construct a Hopf algebra 
on this algebra by adding structure with the  ideas 
coming from the structures of the  affine quantum groups. 
 
The definition of quantum groups  discovered by Drinfeld 
and Jimbo  is presented   as a deformation 
of the  simple Lie algebra
 by the basic generators and the relations based on the data coming from the  
corresponding Cartan matrix. The  extension of the realization of the affine
Kac-Moody algebra $\hat {\frak{g}}$ associated to a simple Lie algebra
$\frak{g}$ as a central extension of the
corresponding loop algebra ${\frak{ g}} \otimes \Bbb C [t, t^{-1}] $
\cite{G} has two different approaches. 
The first approach was given by Faddeev, Reshetikhin and
Takhtajan \cite{FRT}, who obtained a realization  of
the quantum loop algebra $U_q(
{\frak{ g}} \otimes \Bbb [t, t^{-1}])$ via a canonical 
 solution
 of the Yang-Baxter equation depending on a parameter $z\in\Bbb C$. 
This approach was completed by Reshetikhin and Semenov-Tian-Shansky
\cite{RS} by  incorporating  the central extension in the
previous realization. We call this approach FRTS construction. 
The second approach was given  by Drinfeld, who 
\cite{Dr2} gave a realization  
of the quantum affine algebra $U_q(\hgtg)$ and its special degeneration 
called the Yangian. 
As an algebra, this 
realization is equivalent to  the FRTS construction 
 \cite{DF} through
certain Gauss decomposition for the case of $U_q(\hat {\frak gl}(n))$. 
Although we can not   extend the conventional comultiplication to 
the current operators of Drinfeld to derive a closed comultiplication 
formula, 
Drinfeld also gave the Hopf algebra structure 
for such a formulation \cite{DF}, which  
 \cite{DM} \cite{DI}  we used 
to study vertex operators and zeros and poles of the quantum current 
operators. 
In the Drinfeld realization  of 
quantum affine algebras, the structure constants are certain 
rational functions $g_{ij}(z)$. In \cite{DI}, we generalize this type of 
Hopf algebras by substituting  
$g_{ij}(z)$ by other functions that satisfy the functional property 
of $g_{ij}(z)$.   

In this paper, we will  use  the idea of FRTS construction to 
define a Hopf  algebra generated by a current  operator valued matrix on $V$, 
$L(x)$, such that 
$$R(x_1/x_2)L_1(x_1)L_2(x_2)=L_2(x_2)L_1(x_1)R(x_1/x_2), $$
where $L_1(x)=L(x)\otimes 1$ and $L_2(x)=1\otimes L(x)$; 
and the commutation relation between the particles and this
new operator matrix $L(x)$  is presented as 
$$\Phi(x_1)_1 L(x_2)_2=R(x_1/x_2 q^{c/2})^{-1}L(x_2)_2 \Phi(x_1)_1. $$
This relation can be interpreted as that $\Phi(x)$ is an intertwiner for 
the algebra generated by $L(x)$ \cite{FR}. 
With this we can define a comultiplication on the algebra generated by 
$L(x)$ and $\Phi(x)$, where the comultiplication for $L(z)$ comes from 
FRTS construction and the comultiplication of $\Phi(x)$ is 
defined as 
$$\Delta (\Phi(z))= \Phi(x)\otimes 1 +L(x q^{c_1/2})\otimes\Phi(zq^{c_1}),$$ 
which is a generalization of Drinfeld construction. Then combining the 
idea of FRTS construction and Drinfeld realization,  we give a double for 
such a construction, where the FRTS construction is a Hopf subalgebra
and Drinfeld realization is a special case of our realization with 
certain diagonal $R(x)$. This paper is basically a result of the combination 
of the two approaches, the FRTS construction and the Drinfeld realization. 

This paper contains two sections. In the first section, we define the 
Zamolochikov algebra and present its Hopf algebra extension. The second
section is to describe the double of such a construction and the related 
examples.

\section{Hopf algebra extension of a Zamolochikov algebra}

Let V be the vector space $\Bbb C^n$.
Let $x$ be an parameter in $\Bbb C$.  A function valued R-matrix 
 $R(x)$ is an 
function valued operator in 
End($V)\otimes$ End($V)$, which satisfies the so-called Yang-Baxter 
equation: 
$$
R_{12}(z) R_{13}(z/w) R_{23}(w) = R_{23}(w)
R_{13}(z/w) R_{12}(z),
$$
where 
 $R_{12}(x)=f_{ij}(x)\sum_{ij}a_i\otimes b_j\otimes 1=R(x)\otimes 1,
 R_{13}(x)=f_{ij}(x)\sum_{ij} a_i\otimes 1\otimes b_j,  
R_{23}(x)=f_{ij}(x)\sum_{ij}1\otimes a_i\otimes b_j=
1\otimes R,$ and $R(x)=f_{ij}(x)\sum_{ij}a_i\otimes b_j.$
We also require that  $R(x)$ satisfies the unitary condition 
$$ 
R_{21}(z)^{-1} = R({z^{-1}}), $$
where $R_{21}(z)=f_{ij}(x)\sum_{ij}b_j\otimes a_i $. 

\begin{defn}
The associative algebra $P[R(x)]$ is an algebra generated by 
operators $\Phi_i(x)$ index by a linear independent basis
$e_i$ of $V$. Let $\Phi(x)=\Sigma  \Phi_i(x)\otimes e_i$. 
The commutation relations are presented as: 
$$  R(x_1/x_2)\Phi(x_1)_1\Phi(x_2)_2=  
\Phi(x_2)_2\Phi(x_1)_1,  $$
where  $\Phi(x_1)_1\Phi(x_2)_2=\Sigma 
 \Phi_i(x_1)\Phi_j(x_2)e_i\otimes e_j$ and 
$\Phi(x_2)_2\Phi(x_1)_1= \Sigma \Phi_j(x_2)
 \Phi_i(x_1)e_i\otimes e_j$
\end{defn}

As explained in the introduction, this system is used in the  description 
particles ($\Phi_i(x)$) in physics  with the scattering matrix $ R(x)$ 
and in some other context. This relation is also  satisfied by the 
vertex operators for 
quantum affine algebras \cite{FR}  \cite{DFJMN}, etc and 
this type of system also 
appeared in describing the elliptic type of algebras  \cite{FO}. 
However, they are all described as  algebras  not 
 Hopf algebras. Following the idea in \cite{DI}, \cite{B}, 
we would like to extend this algebra with additional current operators
coming from the FRTS construction to give a Hopf algebra structure to 
such a system. 

 \begin{defn}
The  algebra $EP[R(x)]$ is an associative algebra 
generated by $\Phi_i(x)$  indexed  by a linear independent basis
$e_i$ of $V$,  $l_{ij}(x)$ indexed by the linear independent 
basis $e_{ij}$ of End$(V)$ and a central element $k$.  Let 
 $\Phi_(x)= \Phi_i(x)\otimes e_i$ and  the operator valued matrix $L(x)
= \Sigma l_{ij}(x)\otimes e_{ij}$, such that $L(x)$ is invertible.  
 They satisfies the commutation relations:   
$$R(x_1/x_2) \Phi(x_1)_1\Phi(x_2)_2= \Phi(x_2)_2\Phi(x_1)_1, $$
 $$\Phi(x_1)_1L(x_2)_2= R(q^{c/2}  x_1/x_2)^{-1}L(x_2)_2 \Phi_i(x_1)_1, $$
$$ R(x_1/x_2)L(x_1)_1L(x_2)_2= L(x_2)_2 L(x_1)_1 R(x_1/x_2). $$
Here $\Phi(x_1)_1L(x_2)_2=\Sigma  \Phi_i(x_1)L_{kl}(x_2)e_i\otimes e_{kl}$,  
$L(x_2)_2 \Phi_i(x_1)_1= \Sigma L_{kl}(x_2)  \Phi_i(x_1)e_i\otimes e_{kl}$,
 $L(x_1)_1L(x_2)_2= PL(x_1)_2L(x_2)_1P= \Sigma   L_{ij}(x_1)  L_i(x_2)e_{ij}
\otimes e_{kl}, $ and $P$ is the permutation operator. 
\end{defn}

\begin{thm}
The algebra $EP[R(x)]$
has a Hopf algebra structure, which are given 
by the following formulae. 

\noindent{\bf Coproduct $\Delta$}
\begin{align*}
\text{(0)}& \quad \Delta(q^c)=q^c\otimes q^c, \\
\text{(1)}& \quad \Delta(\Phi_i(z))=\Phi_i(z)\otimes 1+
           \Sigma  L_{ij}(zq^{\frac{c_1}{2}})\otimes \Phi_j(zq^{c_1}), \\
\text{(2)}& \quad \Delta(L_{ij}(z))=
            \Sigma L_{ik}(zq^{-\frac{c_2}{2}})\otimes L_{kj}
(zq^{\frac{c_1}{2}}), 
\end{align*}
where $c_1=c\otimes 1$ and $c_2=1\otimes c$.

\noindent{\bf Counit $\vep$}
\begin{align*}
\vep(q^c)=1 & \quad \vep(L_{ij}(z))=\delta_{ij}, \\
            & \quad \vep(\Phi_i^{\pm}(z))=0.
\end{align*}
\noindent{\bf Antipode $\quad a$}
\begin{align*}
\text{(0)}& \quad a(q^c)=q^{-c}, \\
\text{(1)}& \quad a(\Phi_i(z))=\Sigma -(L(zq^{-\frac{c}{2}})^{-1})_{ij}
                               \Phi_j(zq^{-c}), \\
\text{(2)}& \quad a(L(z))=(L(z))^{-1}. 
\end{align*}

\end{thm}

We will use the notation to denote the comultiplication. 
 $$  \Delta \Phi(x_1)= \Phi(x_1)\bar \otimes 1+
            L(x_1q^{\frac{c_1}{2}}) \bar \otimes \Phi(x_1q^{c_1});$$
$$\Delta(L(x_2)= (L(x_2q^{-\frac{c_2}{2}})
\bar \otimes L(x_2q^{\frac{c_1}{2}}))$$

 \noindent{Proof.}
For the comultiplication above we have that 
$$  \Delta \Phi(x_1)_1\Delta L(x_2)_2=$$
$$ (\Phi(x_1)\bar \otimes 1+
            L(x_1q^{\frac{c_1}{2}})\bar \otimes \Phi(x_1q^{c_1}))_1(
             L(x_2q^{-\frac{c_2}{2}}) 
\bar \otimes L(x_2 q^{\frac{c_1}{2}})_2= $$
$$R(x_1/x_2 q^{\frac {c_1+c_2} {2}})^{-1}
  (L(x_2q^{-\frac{c_2}{2}}) \bar \otimes L(x_2q^{\frac{c_1}{2}}))_2 
(\Phi(z)\bar \otimes 1)_1+ $$ 
$$
 R(x_2/x_1 q^{\frac {c_1+c_2} {2}})^{-1} 
(L(x_2q^{-\frac{c_2}{2}})\bar \otimes L(x_2q^{\frac{c_1}{2}}))_2 
(L(x_1q^{\frac{c_1}{2}}) \bar \otimes \Phi(x_1q^{c_1}))_1. $$

$$R(x_1/x_2)\Delta \Phi(x_1)_1 \Delta \Phi(x_2)_2= $$
$$R(x_1/x_2)(\Phi(x_1)\bar \otimes 1+
            L(x_1q^{\frac{c_1}{2}})\bar \otimes \Phi(x_1q^{c_1}))_1
(\Phi(x_2)\bar \otimes 1+
            L(x_2q^{\frac{c_1}{2}}) \bar \otimes \Phi(x_2q^{c_1}))_2= $$
$$(\Phi(x_2)\bar \otimes 1)_2(\Phi(x_1)\bar \otimes 1)_1 + 
(L(x_2q^{\frac{c_1}{2}})\bar \otimes \Phi(x_2q^{c_1}))_2(\Phi(x_1)\bar \otimes 1)_1
+$$ 
$$( L(x_1q^{\frac{c_1}{2}})\bar \otimes \Phi(x_1q^{c_1}))_1 
 ( L(x_2q^{\frac{c_1}{2}})\bar \otimes \Phi(x_2q^{c_1}))_2 + 
 R_{21}(x_2/x_1)^{-1}  L(x_1q^{\frac{c_1}{2}})\bar \otimes \Phi(x_1q^{c_1}))_1
(\Phi(x_2)\bar \otimes 1)_2=$$ 
$$ \Delta \Phi(x_2)_2\Phi(x_1)_1 .$$

This construction of 
comultiplication  follows partially the idea of constructing comultiplications
for the quantum Lie algebra \cite{B}, where the cases 
without the parameter $x$ are 
given. With our construction, we can extend the Hopf algebra structures
to the special  Zamolochikov algebra $Z_{n,k}(\xi, \tau)$, which is defined 
as the algebra generated $\Phi(z)$ with  an Belavin elliptic 
R-matrix  $R(z)$\cite{FO}. We expect that 
the new  Hopf algebra 
structure should be very useful in the study of the 
representation theory  of the elliptic Zamolochikov algebras and hopefully 
even the related Sklyanin elliptic algebras.

\section{ The Double of $EP[R(x)]$}

In this section, we will present a double of the algebra  $EP[R(x)]$ following 
the Drinfeld realization of the quantum affine algebra 
$U_q(\hat {\frak sl}(2))$.

 \begin{defn}
The  algebra $DEP[R(x)]$ is an associative algebra 
generated by $\Phi_i(x)$  indexed  by a linear independent basis
$e_i$ of $V$,  $l_{ij}(x)$ and $l^*_{ij}(x)$ indexed by the linear independent 
basis $e_{ij}$ of End$(V)$, $\Phi^*(x)$ indexed  by a linear independent basis
$e_i^*$ of $V^*$, the dual space of $V$, and a central element $k$.  Let 
 $\Phi(x)= \Phi_i(x)\otimes e_i$, 
 $\Phi^*(x)= \Phi_i^*(x)\otimes e_i^*$ the operator valued matrix $L(x)
= \Sigma l_{ij}(x)\otimes e_{ij}$,$L^*(x)
= \Sigma l^*_{ij}(x)\otimes e_{ij}$,
 such that $L(x)$ and $L^*(x)$ are invertible.  
 They satisfies the commutation relations:   
$$R(x_1/x_2) \Phi(x_1)_1\Phi(x_2)_2= \Phi(x_2)_2\Phi(x_1)_1, $$
 $$\Phi(x_1)_1L(x_2)_2= R(q^{c/2}  x_1/x_2)^{-1}L(x_2)_2 \Phi_i(x_1)_1, $$
$$ R(x_1/x_2)L(x_1)_1L(x_2)_2= L(x_2)_2 L(x_1)_1 R(x_1/x_2). $$
$$ R(x_1/x_2)L^*(x_1)_1L^*(x_2)_2= L^*(x_2)_2 L^*(x_1)_1 R(x_1/x_2), $$
$$R(x_1/x_2q^{-c}) L(x_1)_1L^*(x_2)_2=  L^*(x_2)_2 L^*(x_1)_1 R(x_1/x_2q^c), $$
$$ \Phi^*(x_2)_2\Phi^*(x_1)_1=\Phi^*(x_1)_1\Phi^*(x_2)_2R_{21}(x_2/x_1), $$
$$L^*(x_2)_2\Phi^*(x_1)_1=
\Phi^*(x_1)_1L^*(x_2)_2R_{21}(q^{-c/2}x_2/x_1),$$
$$ \Phi(x_1)_1\Phi(x_2)^*_2-  \Phi^*(x_2)_2\Phi(x_1)_1= 
1/(q-q^{-1})(L^*(wq^{c/2})\delta(z/wq^{-c})- 
\delta(z/wq^c)L^*(zq^{c/2}),  $$
$$L(x_1)_1\Phi^*(x_2)_2 R_{21}(q^{-c/2}x_2/x_1)=\Phi^*(x_2)_2L(x_1)_1, $$
$$R(q^{c/2}
x_1/x_2)L^*(x_1)_1  \Phi(x_2)_2=  \Phi(x_2)_2L^*(x_1)_1. $$  
Here $\Phi(x_1)_1L(x_2)_2=\Sigma  \Phi_i(x_1)L_{kl}(x_2)e_i\otimes e_{kl}$,  
$L(x_2)_2 \Phi_i(x_1)_1= \Sigma L_{kl}(x_2)  \Phi_i(x_1)e_i\otimes e_{kl}$,
 $L(x_1)_1L(x_2)_2= PL(x_1)_2L(x_2)_1P= \Sigma   L_{ij}(x_1)  L_i(x_2)e_{ij}
\otimes e_{kl}, $ and the others are defined in the same way. 
$P$ is the permutation operator.
$\delta(z)$ is the distribution with the support at $1$. 
\end{defn}

\begin{thm}
 $DEP[R(x)]$ has  an Hopf algebra structure. The comultiplication 
$\Delta$, the counit  $\vep$ and  the antipode $\quad a$ are 
given by the following formulas. 

\noindent{\bf Coproduct $\Delta$}
\begin{align*}
\text{(0)}& \quad \Delta(q^c)=q^c\otimes q^c, \\
\text{(1)}& \quad \Delta(\Phi_i(z))=\Phi_i(z)\otimes 1+
           \Sigma  L_{ij}(zq^{\frac{c_1}{2}})\otimes \Phi_j(zq^{c_1}), \\
\text{(2)}& \quad \Delta(L_{ij}(z))=
            \Sigma L_{ik}(zq^{-\frac{c_2}{2}})\otimes L_{kj}
(zq^{\frac{c_1}{2}}), \\
\text{(3)}& \quad \Delta(\Phi^*_i(z))=1\otimes \Phi^*_i(z)+
           \Sigma \Phi^*_j(zq^{c_2})\otimes 
L^*_{ij}(zq^{\frac{c_2}{2}}), \\
\text{(2)}& \quad \Delta(L^*_{ij}(z))=
            \Sigma L^*_{ik}(zq^{\frac{c_2}{2}})\otimes L^*_{kj}
(zq^{-\frac{c_1}{2}}).
\end{align*}
where $c_1=c\otimes 1$ and $c_2=1\otimes c$.

\noindent{\bf Counit $\vep$}
\begin{align*}
\vep(q^c)=1 & \quad \vep(L(z))=\vep(L^*(z))=I, \\
            & \quad \vep(\Phi(z))=0=\vep(\Phi^*(z)) .
\end{align*}

\noindent{\bf Antipode $\quad a$}
\begin{align*}
\text{(0)}& \quad a(q^c)=q^{-c}, \\
\text{(1)}& \quad a(\Phi(z))=-L(zq^{-\frac{c}{2}})^{-1}
                               \Phi(zq^{-c}), \\
\text{(2)}& \quad a(\Phi^*(z))=-\Phi^*(zq^{-c})
                               L^*(zq^{-\frac{c}{2}})^{-1}, \\
\text{(3)}& \quad a(L(z))=L(z)^{-1}, \\
\text{(4)}& \quad a(L^*(z))=L^*(z)^{-1}.
\end{align*}
\end{thm}

{\bf Proof}. 
$$
\Delta \Phi^*(z)_1 \Delta \Phi(w)^*_2 R_{21}(w/z)=$$ 
$$(1\otimes \Phi^*(z)+ \Phi^*(zq^{c_2}) \otimes L^*(zq^{\frac{c_2}{2}}))_1
(1\otimes \Phi^*(w)+ \Phi^*(wq^{c_2}) \otimes L^*(wq^{\frac{c_2}{2}}))_2
R_{21}(x_1/x_2)= $$
$$(1\otimes \Phi^*(w)_2(1\otimes \Phi^*(z)_1+ 
(\Phi^*(wq^{c_2}) \otimes L^*(wq^{\frac{c_2}{2}}))_2
(1\otimes \Phi^*(z))_1+$$ 
$$ ( \Phi^*(zq^{c_2}) \otimes L^*(zq^{\frac{c_2}{2}}))_1
(1\otimes \Phi^*(w))_2 R(z/w)^{-1} +
(1\otimes \Phi^*(w)+ \Phi^*(wq^{c_2}) \otimes L^*(wq^{\frac{c_2}{2}}))_2\times
$$
$$(1\otimes \Phi^*(z)+ \Phi^*(zq^{c_2}) \otimes L^*(zq^{\frac{c_2}{2}}))_1
+ (\Phi^*(wq^{c_2}) \otimes L^*(wq^{\frac{c_2}{2}}))_2
 (\Phi^*(zq^{c_2}) \otimes L^*(zq^{\frac{c_2}{2}}))_1.  $$

$$\Delta \Phi^*(z)_1 \Delta L^*(w)_2 R_{21}(q^{-(c_1+c_2)/2}w/z)= $$
$$(1\otimes \Phi^*(z)+ \Phi^*(zq^{c_2}) \otimes L^*(zq^{\frac{c_2}{2}}))_1
( L^*(wq^{\frac{c_2}{2}})\otimes L^*(wq^{-\frac{c_1}{2}}))_2
R_{21}(q^{-(c_1+c_2)/2}w/z)= $$
$$  \Delta L^*(w)_2 \Delta \Phi^*(z)_1. $$

$$\Delta L(z)_1 \Delta 
\Phi^*(w)_2 R_{21}(q^{-(c_1+c_2)/2}w/z)=$$
$$( L(zq^{-\frac{c_2}{2}})\otimes L
(zq^{\frac{c_1}{2}})_1 (1\otimes \Phi^*(w)+ \Phi^*(wq^{c_2})
 \otimes L^*(wq^{\frac{c_2}{2}}))_2 R_{21}(q^{-(c_1+c_2)/2}w/z)= $$
$$ (1\otimes \Phi^*(w))_2
( L(zq^{-\frac{c_2}{2}})\otimes L
(zq^{\frac{c_1}{2}})_1+ ( L(zq^{-\frac{c_2}{2}})\otimes L
(zq^{\frac{c_1}{2}})_1(\Phi^*(wq^{c_2})
 \otimes L^*(wq^{\frac{c_2}{2}}))_2  R_{21}(q^{-(c_1+c_2)/2}w/z)
= $$ 
$$ (1\otimes \Phi^*(w))_2
( L(zq^{-\frac{c_2}{2}})\otimes L
(zq^{\frac{c_1}{2}})_1+(\Phi^*(wq^{c_2})
 \otimes L^*(wq^{\frac{c_2}{2}}))_2( L(zq^{-\frac{c_2}{2}})\otimes L
(zq^{\frac{c_1}{2}})_1=  
$$ 
$$\Delta (\Phi^*(w))_2\Delta (L(x_1))_1, $$

$$R(q^{(c_1+c_2)/2}z/w)\Delta L^*(z)_1  \Delta \Phi(w)_2=$$
$$R(q^{(c_1+c_2)/2}z/w)
 ( L^*(zq^{\frac{c_2}{2}})\otimes L^*(zq^{-\frac{c_1}{2}}))_1
(\Phi(w)\otimes 1+L(wq^{\frac{c_1}{2}})\otimes \Phi(wq^{c_1}))_2
= $$ 
$$ (\Phi(w)\otimes 1)_2( L^*(zq^{\frac{c_2}{2}})\otimes L^*(zq^{-\frac{c_1}{2}}))_1 + 
R_{21}(q^{-(c_1+c_2)/2}w/z) 
 ( L^*(zq^{\frac{c_2}{2}})\otimes L^*(zq^{-\frac{c_1}{2}}))_1
L(wq^{\frac{c_1}{2}})\otimes \Phi(wq^{c_1}))_2= $$ 
$$  \Delta  \Phi(w)_2 \Delta L^*(z)_1. $$

$$\Delta \Phi(z)_1 \Delta \Phi(w)^*_2-\Delta  \Phi^*(w)_2
\Delta \Phi(z)_1= $$
$$(\Phi(z)\otimes 1+L(zq^{\frac{c_1}{2}})\otimes \Phi(zq^{c_1}))_1
(1\otimes \Phi^*(w)+ \Phi^*(wq^{c_2}) \otimes L^*(wq^{\frac{c_2}{2}}))_2 - $$
$$(1\otimes \Phi^*(w)+ \Phi^*(wq^{c_2}) \otimes L^*(wq^{\frac{c_2}{2}}))_2
(\Phi(z)\otimes 1+L(zq^{\frac{c_1}{2}})\otimes \Phi(zq^{c_1}))_1=$$
$$ 0+ 1/(q-q^{-1})(L^*(wq^{c_1/2+c_2})\delta(z/wq^{-c_1-c_2})\otimes 
 L^*(wq^{\frac{c_2}{2}})- 
\delta(z/wq^{c_1-c_2})L(zq^{c_1/2}))\otimes L^*(wq^{\frac{c_2}{2}})+ $$ 
$$
 L(zq^{\frac{c_1}{2}})\otimes(1/(q-q^{-1})(L^*(wq^{\frac{c_2}{2}})
\delta(z/wq^{-c_2
+c_1})- 
\delta(z/wq^{c_1+c_2})L(zq^{c_1+c_2/2})+ $$
$$(L(zq^{\frac{c_1}{2}})\otimes \Phi(zq^{c_1}))_1 (
 \Phi^*(wq^{c_2}) \otimes L^*(wq^{\frac{c_2}{2}}))_2- 
  ( \Phi^*(wq^{c_2}) \otimes L^*(wq^{\frac{c_2}{2}}))_2
(L(zq^{\frac{c_1}{2}})\otimes \Phi(zq^{c_1}))_1$$

Because 

$$(L(zq^{\frac{c_1}{2}}\otimes 1)_1
( \Phi^*(wq^{c_2}))\otimes 1)_2  R_{21}(q^{-c_1+c_2}w/z) =
( \Phi^*(wq^{c_2}))\otimes 1)_2(L(zq^{\frac{c_1}{2}}\otimes 1)_1, $$ 
and 
$$ ( R_{21}(q^{-c_1+c_2}w/z))^{-1}
(1\otimes \Phi(zq^{c_1}))_1 L^*(wq^{\frac{c_2}{2}}))_2 
=  L^*(wq^{\frac{c_2}{2}}))_2(1\otimes \Phi(zq^{c_1}))_1, $$

we have that 

$$\Delta \Phi(z)_1 \Delta \Phi(w)^*_2-\Delta  \Phi^*(w)_2
\Delta \Phi(z)_1=$$ 
$$
1/(q-1^{-1})(\Delta L(z/wq^{(c_1+c_2)/2}\delta(z/wq^{-(c_1+c_2)})- 
\delta(z/wq^c)\Delta L^*(zq^{(c_1+c_2)/2}).  $$

In all the setting above, we assume 
$l_{ij}(z)$, $l_{ij}^*(z)$, $\Phi_i(z)$ and $\Phi^*_i(z)$
are functional operators, namely the operator depending 
the variable $z$. On the other hand, we can 
assume that $z$ is a formal variable and 
$l_{ij}(z)= \Sigma_{n\in \Bbb Z} l_{ij}(n)z^{-n},$
 $l^*_{ij}(z)= \Sigma_{n\in \Bbb Z} l^*_{ij}(n)z^{-n},$
$\Phi_{i}(z)= \Sigma_{n\in \Bbb Z} \Phi_{i}(n)z^{-n},$
$\Phi^*_{i}(z)= \Sigma_{n\in \Bbb Z} \Phi^*_{i}(n)z^{-n}.$
We can define an algebra  $DZP[R(x)]$.

Let $R'(z)= R(z)f(z)$, where $f(z)$ is the common divisor of all the functions
$F(z)$, such that $F(z)  R(z)$ has no poles.
\begin{defn}
The  algebra $DZP[R(x)]$ is an associative algebra 
generated by $\Phi_i(x)$  indexed  by a linear independent basis
$e_i$ of $V$,  $l_{ij}(x)$ and $l^*_{ij}(x)$ indexed by the linear independent 
basis $e_{ij}$ of End$(V)$, $\Phi^*(x)$ indexed  by a linear independent basis
$e_i^*$ of $V^*$, the dual space of $V$, and a central element $k$.  Let 
 $\Phi(x)= \Phi_i(x)\otimes e_i$, 
 $\Phi^*(x)= \Phi_i^*(x)\otimes e_i^*$ the operator valued matrix $L(x)
= \Sigma l_{ij}(x)\otimes e_{ij}$,$L^*(x)
= \Sigma l^*_{ij}(x)\otimes e_{ij}$,
 such that $L(x)$ and ,$L^*(x)$ are invertible.
 They satisfies the commutation   
$$R'(x_1/x_2) \Phi(x_1)_1\Phi(x_2)_2=f(x_1/x_2) \Phi(x_2)_2\Phi(x_1)_1, $$
 $$\Phi(x_1)_1L(x_2)_2= R(q^{c/2}  x_1/x_2)^{-1}L(x_2)_2 \Phi_i(x_1)_1, $$
$$ R(x_1/x_2)L(x_1)_1L(x_2)_2= L(x_2)_2 L(x_1)_1 R(x_1/x_2). $$
$$ R(x_1/x_2)L^*(x_1)_1L^*(x_2)_2= L^*(x_2)_2 L^*(x_1)_1 R(x_1/x_2), $$
$$R(x_1/x_2q^{-c}) L(x_1)_1L^*(x_2)_2=  L^*(x_2)_2 L^*(x_1)_1 R(x_1/x_2q^c), $$
$$f(x_2/x_1)\Phi^*(x_2)_2\Phi^*(x_1)_1=\Phi^*(x_1)_1\Phi^*(x_2)_2R'_{21}(x_2/x_1), $$
$$L^*(x_2)_2\Phi^*(x_1)_1=
\Phi^*(x_1)_1L^*(x_2)_2R_{21}(q^{-c/2}x_2/x_1),$$
$$ \Phi(x_1)_1\Phi(x_2)^*_2-  \Phi^*(x_2)_2\Phi(x_1)_1= 
1/(q-q^{-1})(L^*(wq^{c/2})\delta(z/wq^{-c})- 
\delta(z/wq^c)L^*(zq^{c/2}),  $$
$$L(x_1)_1\Phi^*(x_2)_2 R_{21}(q^{-c/2}x_2/x_1)=\Phi^*(x_2)_2L(x_1)_1, $$
$$R(q^{c/2}
x_1/x_2)L^*(x_1)_1  \Phi(x_2)_2=  \Phi(x_2)_2L^*(x_1)_1. $$  
Here $\Phi(x_1)_1L(x_2)_2=\Sigma  \Phi_i(x_1)L_{kl}(x_2)e_i\otimes e_{kl}$,  
$L(x_2)_2 \Phi_i(x_1)_1= \Sigma L_{kl}(x_2)  \Phi_i(x_1)e_i\otimes e_{kl}$,
 $L(x_1)_1L(x_2)_2= PL(x_1)_2L(x_2)_1P= \Sigma   L_{ij}(x_1)  L_i(x_2)e_{ij}
\otimes e_{kl}, $ and the others are defined in the same way as above. 
$P$ is the permutation operator.
$\delta(z)=\Sigma_{n\in \Bbb Z} z^n$. The operator $R(z)$ and $R_21(z)$ are 
expanded in appropriate directions. 
\end{defn}

If the poles of the matrix of $R(z)$ are beyond a finite disc around zero, 
we can always impose the condition that $l_{kl}(n)=0=l^*_{kl}(-n)=
l_{ij}(0)=l^*_{ji}(0)$, for $n<0$, $i<j$. Then the condition of the 
invertibility is equivalent to requires that $l_{ii}(0)$ and 
$l_{ii}(0)$ are invertible.

\begin{ex}
Let $v$ be one dimensional, and  $R(z)= z-wq^2/zq^2-w$. 
Let $l_{11}(n)=0=l^*_{11}(-n)$, for $n<0$. 
Then the algebra  $DZP[R(x)]$ is the quantum affine algebra 
$U_q(\hat {\frak sl}(2))$. 
If we choose   $R(z)$ to be other functions with the property 
$R(z)= (R(z^{-1})^{-1}$, then it is an algebra defined in 
\cite{DI} as an generalization of the 
 the quantum affine algebra 
$U_q(\hat {\frak sl}(2))$. 
\end{ex}

\begin{ex}
Let $v=C^n$ and  $R(z)=\Sigma 
 (z-wq^2)/(zq^2-w)e_{ii}\otimes e_{ii}+ \Sigma 
 (z-wq^{-1})
/(zq^{-1}-w)(e_{ii}\otimes e_{i+1,i+1}+ e_{i+1,i+1} \otimes e_{i,i}).$ 
Let $l_{kl}(n)=0=l^*_{kl}(-n)=l_{ij}(0)=l^*_{ji}(0)$, for $n<0$, $i<j$. 
Then the algebra  $DZP[R(x)]$ is an algebra, whose quotient 
( modular the cubic relations) is  
$U_q(\hat {\frak sl}(n))$. 
If we substitute  $ z-wq^2/zq^2-w$ and  $(z-wq^{-1})
/(zq^{-1}-w)$ by other functions, it will be the generalization of 
$U_q(\hat {\frak sl}(n))$ without  the cubic relations \cite{DI}. 
\end{ex}

\begin{ex}
Let $v=C^n$ and  $R(z)$ be the 
projection of the universal R-matrix ${\frak R} $
of $U_q(\hat {\frak sl}(n))$ . 
Let $l_{kl}(n)=0=l^*_{kl}(-n)=l_{ij}(0)=l^*_{ji}(0)$, for $n<0$, $i<j$. The
operator $L(z)$ can be identified with the operator 
$(id\otimes \pi_V) {\frak R}_21(zq^{c/2})$ and
$L^*(z)$  with the operator 
$(id\otimes \pi_V) {\frak R}^{-1}(z^{-1}q^{-c/2})$. 
The subalgebra generated by $L(z)$ and $L^*(z)$ is isomorphic to 
$U_q(\hat {\frak gl}(n))$.  It can  see that the algebra 
 $DZP[R(x)]$ should be isomorphic to $U_q(\hat {\frak sl}(n+1))$, 
because when $q$ goes to 1, this algebra degenerate into 
 $\hat {\frak sl}(n+1)$. 
\end{ex} 

From the definition, we can see both  the subalgebra generated by 
$\Phi(z)$, $L(z)$ and $L^*(z)$ and the subalgebra generated by  
$\Phi^*(z)$, $L(z)$ and $L^*(z)$ are the Hopf algebras. If we take 
$R(z)$ to be the projection of the universal R-matrix ${\frak R} $
of $U_q(\hat {\frak sl}(n))$ on certain linear spaces, we will derive 
unconventional Hopf algebras from those subalgebras. 

It is clear our new algebras can be viewed as a 
 simple generalization of the 
Drinfeld realization of the quantum affine algebra
$U_q(\hat {\frak sl}(2))$, where the function 
$g(z)$ is substitute by a matrix $R(z)$,  the operators 
are substituted  by  the vector valued operators and the relations 
looks the same. However such a generalization is highly non-trivial in the 
sense that all the Hopf algebra structures are preserved, in the other 
words, those new algebras are Hopf algebras, whose comultiplication, 
counit and antipode symbolically are the same.  These new Hopf algebras 
should be very useful  in various applications  in mathematics and 
physics, for example,  the study of the 
representation theory  of the elliptic Zamolochikov algebras.

\begin{ack}
We would like to thank M. Jimbo 
and B. Feigin for useful discussions. 
This project 
is supported by the grant Reward research (A) 08740020 from the 
Ministry of Education of Japan. 
\end{ack}


\begin{thebibliography}{[Beck]}

\bibitem [B] {B}
D. Bernard 
{\it Quantum Lie algebra and differential calculus on quantum groups},
Preprint, SPht-90-119

\bibitem [DFJMN] {DFJMN}
{\it Diagonalization of the XXZ Hamiltonian by vertex operators}, 
CMP, {\bf 151}, 1993, 89-153


\bibitem   [DF] {DF} J. Ding, I. B. Frenkel
{\it Isomorphism of two realizations of quantum affine algebra 
 $U_q(\hat {\frak gl}(n))$ },
 CMP, {\bf 156}, 1993, 277-300
Physics 

\bibitem   [DI] {DI} J. Ding, K. Iohara
{\it Generalization and deformation of 
the quantum affine algebras}, 
 RIMS-1090

\bibitem [DM] {DM} J.Ding and T. Miwa 
{\it Zeros and poles of quantum current operators and 
the quantum integrable condition}, RIMS-1092. 

\bibitem    [Dr1] {Dr1} V. G. Drinfeld 
{\it Hopf algebra and the quantum Yang-Baxter Equation},
Dokl. Akad. Nauk. SSSR, {\bf 283}, 1985, 1060-1064

 
\bibitem  [Dr2]{Dr2}  V. G. Drinfeld {\it
 New realization of Yangian and quantum
affine algebra}, Soviet Math. Doklady,
{\bf 36}, 1988,  212-216

 \bibitem [FRT] {FRT}
 L. D. Faddeev, N. Yu, Reshetikhin, L. A. Takhtajan
{\it
Quantization of Lie groups and Lie algebras,  Yang-Baxter equation in
Integrable Systems},   (Advanced Series in Mathematical Physics {\bf 10})
World Scientific, 1989, 299-309.

\bibitem [FO] {FO}  B. Feigin, V. Odesski
{\it  Sklyanin Elliptic algebras}
Funkts. Anal. Prilozhen., {\bf 23}, 3, 1989, 45-54

\bibitem [FO1] {FO1}  B. Feigin, V. Odesski
{\it Vector bundles on Elliptic curve and Sklyanin algebras}
RIMS-1032, q-alg/9509021

\bibitem [FF] {FF} E. Frenkel, B. Feigin
{\it Quantum ${\cal W}$-algebra and elliptic algebras}
RIMS-1027, q-alg/9508009.

\bibitem [FR] {FR} I. B.  Frenkel, N. Yu, Reshetikhin
{\it Quantum affine algebras and holonomic difference equations}, 
CMP, {\bf 146}, 1992, 1-60 

\bibitem [G] {G}  H. Garland {\it
  The arithmetic theory of loop groups},
Publ. Math. IHES {\bf 52}, 1980,  5-136


\bibitem     [J1]{J1} M. Jimbo 
{\it A $q$-difference analogue of $U({\frak g})$ and Yang-Baxter equation
},
Lett. Math. Phys. {\bf 10}, 1985, 63-69

\bibitem
[RS] {RS}  N.Yu. Reshetikhin, M.A. Semenov-Tian-Shansky
{\it  Central Extensions of Quantum Current Groups}, LMP, 
{\bf 19}, 1990

\bibitem
[S] {S}  E. K. Sklyanin {\it 
On some algebraic structures related to the 
Yang-Baxter equation} 
Funkts. Anal. Prilozhen, {\bf 16}, No. 4, 1982, 22-34


\end{thebibliography}
\end{document}